\documentclass{article}
\usepackage{amsfonts}

\def\rest{\mathord{\restriction}}
\def\BbbN{\mathbb N}
\def\qed{\hspace*{0pt}\hfill
\hbox{\begin{picture}(10,10)\put(0,0){\framebox(6,6){}}
\end{picture}}\medskip}

\renewcommand{\phi}{\varphi}

\newcommand{\su}{\subseteq}
\renewcommand{\a}{\alpha}
\renewcommand{\b}{\beta}

\newcommand{\om}{\omega}

\newcommand{\lng}{\langle}
\newcommand{\rng}{\rangle}

\newcommand{\sm}{\setminus}

\newcommand{\Cal}[1]{{\mathcal{#1}}}
 \usepackage{amstext}
 \usepackage{amssymb}

\newtheorem{lemma}{Lemma}
\newtheorem{theorem}[lemma]{Theorem}

\renewcommand{\int}{\operatorname{\rm int}}

\newenvironment{proof}{\noindent{\bf Proof}.\ }{\qed}

\begin{document}

\title{van der Waerden spaces and Hindman spaces are not the same}

\author{Menachem Kojman\\
Department of
Mathematics\\
Ben Gurion University of the Negev\\
Beer Sheva, Israel\\
E-mail: {\tt kojman@cs.bgu.ac.il}
\thanks{The first author was partially supported by an Israel
Science Foundation grant.}
\and Saharon Shelah\\
Institute of Mathematics\\
The Hebrew University of Jerusalem\\
Jerusalem, Israel\\
E-mail: {\tt shelah@ma.huji.ac.il}
\thanks{The second author was partially supported by an Israel
Science Foundation grant. Number 781 in Shelah's list of publications}
\thanks{The authors wish to acknowledge a substantial simplification
made by the referee in the proof. The referee has eliminated an
inessential use that the authors have made of the canonical van der
Waerden theorem, all of whose known proofs use Szemer\'edi's theorem.}}

\maketitle

\begin{abstract}
A Hausdorff topological space $X$ is \emph{van der Waerden} if
for every sequence $(x_n)_{n\in\omega}$ in $X$ there is a converging
subsequence $(x_n)_{n\in A}$ where $A\su \om$ contains arithmetic
progressions of all finite lengths. A Hausdorff topological space
$X$ is \emph{Hindman} if for every sequence $(x_n)_{n\in\omega}$ in $X$
there is an \emph{IP-converging} subsequence $(x_n)_{n\in FS(B)}$
for some infinite $B\su \om$.

We show that the continuum hypothesis implies the existence of a
van der Waerden space which is not Hindman.
\end{abstract}

\section{Introduction}

A Hausdorff topological space $X$ is \emph{van der Waerden} if
for every sequence $(x_n)_{n\in\omega}$ in $X$ there is a converging
subsequence $(x_n)_{n\in A}$ where $A\su \om$ contains arithmetic
progressions of all finite lengths. A Hausdorff topological space
$X$ is \emph{Hindman} if for every sequence $(x_n)_{n\in\omega}$ in $X$
there is an \emph{IP-converging} subsequence $(x_n)_{n\in FS(B)}$
for some infinite $B\su \om$. The term $FS(B)$ stands for the set
of all \emph{finite sums} (with no repetitions) over $B$ and
IP-convergence to a point $x\in X$ means: for every neighborhood
$U$ of $x$,
there is some $n_0$ so that $\{x_n:n\in
FS(B\sm \{0,1,\ldots,n_0-1\})\}\su U$.

The classes of van der Waerden  and of Hindman spaces were
introduced in \cite{K1,K2} where it was shown that each class was
productive and properly contained in the class of sequentially
compact spaces, and  that every Hausdorff space $X$ in which the
closure of every countable set is compact and first countable is
both van der Waerden and Hindman. The question was raised whether
every Hausdorff space $X$ is van der Waerden if and only if it is
Hindman. We answer this question in the negative using the
Continuum Hypothesis.

\subsection{Notation and combinatorial preliminaries}

A set $A\su \om$ is an \emph{AP-set} if it contains arithmetic
progressions of all finite lengths. By van der Waerden's theorem
\cite{VDW}, if an AP-set $A$ is partitioned into finitely many
parts, at least one of the parts is AP. Let $\mathcal I_{AP}$
denote the collection of all subsets of $\om$ which are not AP.
$\mathcal I_{AP}$ is a proper ideal over $\om$ and a set $A\su
\om$ is AP if and only if $A\notin \mathcal I_{AP}$.

A set $A\su \om$ is an \emph{IP-set} if there exists an infinite
set $B\su \om$ so that $FS(B)\su A$. $FS(B)=\{\sum F: F\su A,
|F|<\aleph_0\}$, where $\sum F$ stands for $\sum_{n\in F} n$. By
Hindman's theorem \cite{Hindman}, if an IP-set  $A$ is
partitioned into finitely many parts, at least one of the parts
is IP. Let $\mathcal I_{IP}$ denote the collection of all subsets
of $\om$ which are not IP. $\mathcal I_{IP}$ is a proper ideal
over $\om$ and a set $A\su \om$ is IP if and only if $A\notin
\mathcal I_{IP}$.

 We shall need the following lemma which
relates $\mathcal I_{AP}$  to $\mathcal I_{IP}$.

\begin{lemma} \label{apnip}
Let $A$ be an AP set and let $f:\omega\to\omega$.
There exists an AP set $C\subseteq A$ such that either
\begin{enumerate}
\item[(1)] $|f[C]|=1$ or
\item[(2)] $f$ is finite-to-one on $C$
and if $\langle x_n\rangle_{n=0}^\infty$ enumerates $f[C]$ in
increasing order, then $\displaystyle \lim_{n\to\infty}(x_{n+1}-x_n)
=\infty$.
\end{enumerate}
\noindent
In particular, $f[C]\in{\cal I}_{IP}$.
\end{lemma}

\begin{proof}
Suppose that for every AP set $C\subseteq A$,
$|f[C]|>1$.  We construct an AP set $C\subseteq A$ for which
conclusion (2) holds.

For each $m\in\omega$, $A\cap f^{-1}[\{0,1,\ldots,m-1\}]$ is
not an AP set because it is the finite union of sets on which
$f$ is constant, and thus $A\setminus f^{-1}[\{0,1,\ldots,m-1\}]$
is an AP set.  (We are using here the fact that when an AP set
is partitioned into finitely many parts, one of these parts is
an AP set.)

We construct inductively sets $C_n$ for each $n\in\BbbN$ such that
\begin{enumerate}
\item[(a)] for each $n\in\BbbN$, $C_n$ is a length $n$ arithmetic
progression and
\item[(b)] for all $n,m\in\BbbN$, all $x\in C_m$, and all
$y\in C_n$,
if $m< n$, then $f(y)\geq f(x) + n$ and
if $m=n$, then either $f(x)=f(y)$ or
$|f(x)-f(y)|\geq n$.
\end{enumerate}
Let $C_1$ be any singleton subset of $A$.  Let $n\in\BbbN$ and
assume that we have chosen $C_1,C_2,\ldots,C_n$.
Let $k=\max\bigcup_{i=1}^n f[C_i]$ and
choose $i\in\{0,1,\ldots,n\}$ such that\break
$(A\setminus f^{-1}[\{0,1,\ldots,k+n\}])\cap f^{-1}[(n+1)\omega+i]$ is
an AP set. Let $C_{n+1}$ be a length $n+1$ arithmetic
progression contained in
$(A\setminus f^{-1}[\{0,1,\ldots,k+n\}])\cap f^{-1}[(n+1)\omega+i]$.
Given $m\leq n+1$, $x\in C_m$, and $y\in C_{n+1}$, if
$m\leq n$, then $f(x)\leq k$ and $f(y)\geq k+n+1$, while if
$m=n+1$, then either $f(x)=f(y)$ or $|f(x)-f(y)|\geq n+1$.

Let $C=\bigcup_{n=1}^\infty C_n$.
\end{proof}

\section{The space}

\begin{lemma}\label{mad}
Assume CH. Then there exists a maximal  almost disjoint
 family $\mathcal A\su
\mathcal I_{IP}$ so that for every AP-set $B\su \om$ and every
finite-to-one function $f:B\to \om$ there exists an AP-set $C\su
B$ and $A\in \mathcal A$ so that $f[C]\su A$.
\end{lemma}

\begin{proof}
We construct from CH an
 almost disjoint family $\Cal A = \{A_\a:\a<\om_1\}\su
\mathcal I_{IP}$ by induction on $\a$. The enumeration
$\{A_\a:\a<\om_1\}$ may contain repetitions. Let $\{A_n:n<\om\}\su
\mathcal I_{IP}$ be a collection of infinite and pairwise
disjoint sets.

Fix a list $\lng (f_\a,B_\a):\om\le\a<\om_1\rng$ of all pairs
$(f,B)$ in which $B\su \om$ is an AP-set  and $f:B\to \om$ is a
finite-to-one function.

Suppose $\om\le\a<\om_1$ and that $A_\b$ has been chosen for all
$\b<\a$. Consider the pair $(f_\a,B_\a)$.
 If there exists a  finite
set $\{\b_0,\b_1,\dots,\b_\ell\}\su \a$ so that
$f_\a^{-1}[\bigcup_{i\le\ell} A_{\b_i}]$ is AP, let $A_\a=A_0$.

Otherwise, enumerate $\a$ as $\lng \b_i:i<\om\rng$, and now for
all $n<\om$ the set $f_\a^{-1}[\bigcup_{i<n}A_{\b_i}]$ is not AP,
hence $B_\a\sm f_\a^{-1}[\bigcup_{i<n}A_{\b_i}]$ is AP. Let an
arithmetic progression  $D_n\su B_\a\sm
f_\a^{-1}[\bigcup_{i<n}A_{\b_i}]$ of length $n$  be chosen for
all $n$. Then $B':=\bigcup_{n\in\om} D_n$ is  an AP-subset of $B_\a$,
$f_\a[B']$ is infinite (because $f_\a$ is finite-to-one) and
$|f_\a[B']\cap A_\b|<\aleph_0$  for all $\b<\a$. By Lemma
\ref{apnip}  find an AP-set  $B''\su B'$, so that $f_\alpha[B'']\in
\mathcal I_{IP}$, and define $A_\a=f_\alpha[B'']$.

The family $\mathcal A=\{A_\a:\a<\om_1\}$ is clearly an almost
disjoint family of (infinite) sets, and $\mathcal A\su \mathcal
I_{IP}$.

Suppose now that $B\su \om$ is an AP-set and that $f:B\to \om$ is
finite-to-one. There is an index $\om\le\a<\om_1$ for which
$(B,f)=(B_\a,f_\a)$. At stage $\a$ of the construction of
$\mathcal A$, either $f^{-1}[A_{\b_0}\cup\dots\cup A_{\b_\ell}]$
was AP for some finite set $\{\b_0,\dots,\b_\ell\}\su \a$, hence
$f^{-1}[A_\b]$ was AP for some single $\b<\a$, or else
$f^{-1}[A_\a]$ was AP. In either case, there is an AP-set $C\su
B$ and $A\in \mathcal A$ so that $f[C]\su A$.

Finally, to verify that $\mathcal A$ is  maximal let an infinite
set $D\su \om$ be given and let $f:\om\to D$ be the increasing
enumeration of $D$. Since there is an AP-set $C\su \om$ and $A\in
\mathcal A$ so that $f[C]\su A$ it is clear that $D\cap A$ is
infinite.
\end{proof}

\medskip
\begin{theorem} Suppose CH holds. Then there exists a compact,
separable
 van der Waerden space  which is not Hindman.
\end{theorem}

\begin{proof}
Let $\mathcal A$ be as stated in Lemma \ref{mad}. For each $A\in
\mathcal A$ let $p_A\notin \om$ be a distinct point. Define a
topology $\tau$ on $Y=\om\cup\{p_A:A\in \mathcal A\}$ by
requiring that $Z\in \tau$  if and only if for all $p_A\in Z$ the
set $A\sm Z$ is finite. Then for each $A\in{\cal A}$,
$A\cup\{p_A\}$ is a compact neighborhood of $p_A$,
so $\tau$ is a locally compact Hausdorff
topology in which $\om$ is a dense and discrete subspace. Let
$X=Y\cup\{p\}$  be the one-point compactification of $\tau$.

It was shown in \cite[Theorem 10]{K2} that when $\mathcal A\su
\mathcal I_{IP}$ is maximal almost disjoint, the space
constructed in this way is sequentially compact but not Hindman.
To keep this paper self contained, we repeat the simple argument
showing that $X$ is not Hindman.  For each $n\in\om$, let
$x_n=n$ and suppose we have some infinite $B\su \om$ such that
$(x_n)_{n\in FS(B)}$ IP-converges to $q\in X$.  Then $q\notin\om$.
If $q=p_A$ for some $A\in{\cal A}$, then $A$ is an IP set.  So
$q=p$.  By the maximality of ${\cal A}$, pick $A\in{\cal A}$ such
that $A\cap B$ is infinite.  But then $X\sm (A\cup\{p_A\})$ is
a neighborhood of $p$ and for no $n$ does one have
$FS(B\sm\{0,1,\ldots,n-1\})\subseteq X\sm (A\cup\{p_A\})$.

We have yet to see that $X$ is van der Waerden. Suppose $f:\om\to
X$ is given. Let $g:f[\om]\to\om$ be 1-1. By Lemma
\ref{apnip} we can find an AP set $B\su \om$ so that $(g\circ
f)\rest B$ is constant or finite-to-one, and hence $f\rest B$ is
constant or finite-to-one. In the former case, the sequence
$(f(n))_{n\in B}$ is constant, and therefore converges. So assume
that $f\rest B$ is finite-to-one. Since either $f^{-1}[\om]\cap
B$ or $B\sm f^{-1}[\om]$ is AP, we may assume, by shrinking $B$
to some AP-subset, that either $f[B]\su \om$ or $f[B]\su X\sm
(\om\cup\{p\})$.

In the former case,
there is some $A\in \mathcal A$ and AP-set $C\su B$ so that
$f[C]\su A$. Since $f\rest B$ is finite-to-one, $(f(n))_{n\in C}$
converges to $p_A$.
In the latter case, we claim that the
sequence $(f(n))_{n\in B}$ converges to $p$.
To see this, let $Z$ be a compact subset of $Y$, so that
$X\sm Z$ is a basic neighborhood of $p$.  Then $Z\sm\om$ is finite so,
since $f\rest B$ is finite-to-one, $(f(n))_{n\in B}$ is eventually
in $X\sm Z$.
\end{proof}

\end{document}